\UseRawInputEncoding

\documentclass[12pt]{amsart}
\usepackage{amsfonts}
\usepackage{amsmath}
\usepackage{amssymb,latexsym}
\usepackage[mathcal]{eucal}
\usepackage{amscd}

\usepackage[pdftex,bookmarks,colorlinks,breaklinks]{hyperref}
\input xy
\xyoption{all}

\newcommand{\Acal}{\mathcal{A}}
\newcommand{\Bcal}{\mathcal{B}}

\newcommand{\Dcal}{\mathcal{D}}

\newcommand{\Fcal}{\mathcal{F}}

\newcommand{\Lcal}{\mathcal{L}}
\newcommand{\Mcal}{\mathcal{M}}

\newcommand{\Scal}{\mathcal{S}}

\newcommand{\Wcal}{\mathcal{W}}

\newcommand{\Z}{\mathbb{Z}}
\newcommand{\Q}{\mathbb{Q}}
\newcommand{\R}{\mathbb{R}}
\newcommand{\C}{\mathbb{C}}
\newcommand{\N}{\mathbb{N}}
\newcommand{\T}{\mathbb{T}}

\newcommand{\al}{\alpha}
\newcommand{\Ga}{\Gamma}
\newcommand{\ga}{\gamma}
\newcommand{\del}{\delta}
\newcommand{\Del}{\Delta}

\newcommand{\sig}{\sigma}

\newcommand{\la}{\lambda}

\newcommand{\Om}{\Omega}

\newcommand{\br}{\vspace{3 mm}}

\newcommand{\nor}{\vartriangleleft}

\newcommand{\ben}{\begin{enumerate}}
\newcommand{\een}{\end{enumerate}}

\newcommand{\vertiii}[1]{{\left\vert\kern-0.25ex\left\vert\kern-0.25ex\left\vert #1 
    \right\vert\kern-0.25ex\right\vert\kern-0.25ex\right\vert}}

\swapnumbers
\theoremstyle{plain}
\newtheorem{thm}{Theorem}[section]
\newtheorem{cor}[thm]{Corollary}
\newtheorem{lem}[thm]{Lemma}

\theoremstyle{definition}

\newtheorem{rmk}[thm]{Remark}

\newtheorem{claim}[thm]{Claim}

\begin{document}

\title[A family of distal functions]
{A family of distal functions and multipliers for strict ergodicity}

\date{June 20, 2021}


\author{Eli Glasner}

\address{Department of Mathematics\\
     Tel Aviv University\\
         Tel Aviv\\
         Israel}
\email{glasner@math.tau.ac.il}
\begin{abstract}
We give two proofs to an old result of E. Salehi, 
showing that the Weyl subalgebra $\Wcal$ of $\ell^\infty(\Z)$ is a proper subalgebra of $\Dcal$,
the algebra of distal functions. 
We also show that the family 
$\Scal^d$ of strictly ergodic functions in $\Dcal$ does not form an algebra and hence in
particular does not coincide with $\Wcal$.
We then use similar constructions to show that 
a  function which is a multiplier for strict ergodicity,
either within $\Dcal$ or in general, is necessarily a constant.
An example of a metric, strictly ergodic, distal flow 
is constructed which admits a non-strictly ergodic $2$-fold minimal self-joining.
It then follows that the enveloping group of this flow is not strictly ergodic (as a $T$-flow).
Finally we show that the distal, strictly ergodic Heisenberg nil-flow is relatively
disjoint over its largest equicontinuous factor from $|\Wcal|$.
\end{abstract}




%
%
%

\maketitle


\section{Introduction}

This work was originally written in the late 1980s.
It was then circulated as a preprint among a few colleagues (see e.g. MathSciNet review: 986700) but
for some reason was never published. In a recent conversation Hillel Furstenberg asked
whether there exists a  subalgebra of $\ell^\infty(\Z)$ which 
is contained in the subalgebra of distal functions,  consists of strictly ergodic functions,
and is ``universal" in some sense.
In a way Section \ref{sec-se} of this work gives a negative answer to Hillel's question.
So I decided to resurrect this old work.
The new version differs from the original one mostly in the addition of a proof of Theorem 
\ref{ab}, and the new Section \ref{sec-new}. 
Also  some misprints are corrected and a few references are added;
in particular a reference to Salehi's work \cite{S}.


\br

A {\em flow} in this work is a pair $(X,T)$ consisting of a compact space $X$ and a self homeomorphism 
$T : X \to X$. A flow $(X,T)$ is {\em minimal} if the {\em orbit} 
$\{T^n x : n \in \Z\}$ of every point $x \in X$, is dense in $X$.
If $(X,T)$ and $(Y,T)$ are two flows then a {\em homomorphism}
$\pi : (X,T) \to (Y,T)$ is a continuous, surjective map $\pi : X \to Y$ such that
$\pi(Tx) = T\pi(x)$ for every $x \in X$.
(We often use the same letter $T$ to denote the acting homeomorphism 
on various flows; the idea is that these are all $\Z$-actions 
where $T$ represents the generator $1$ of $\Z$.)
The {\em enveloping semigroup} $E(X,T)$ of a flow $(X,T)$ is the
compact sub-semigroup of the {\em right topological semigroup} $X^X$,
formed as the closure of the image of $\Z$ in $X^X$ under the map
$n \mapsto T^n$. $E(X,T)$ is both a compact right topological semigroup and a flow,
were $T$ acts by multiplication on the left.
For more details see e.g.  \cite{Gl-07}.

A metric flow $(X,T)$ is called {\em distal} if 
$x, y \in X, \ x \neq y$, implies $\inf_{n \in \Z} d(T^n x, T^n y) >0$.
(When $X$ is nonmetrizable the requirement is: 
$\forall \ x, y \in X, \ x \neq y$, there is a neighborhood $V$ of the diagonal $\Del \subset X \times X$
such that $(T^nx, T^ny) \not \in V, \ \forall n \in \Z$.)
By basic theorems of Robert Ellis \cite{Ellis-58} (i) Every distal flow $(X,T)$ is 
 {\em semisimple}
(i.e. it is the union of its minimal subsets),
and  (ii) A flow $(X,T)$ is distal iff its enveloping semigroup $E(X,T)$ is a group.

A theorem of Krylov and Bogolubov asserts that every $\Z$-flow admits an invariant,  Borel,
probability measure (see e.g. \cite[Theorem 4.1]{Gl2}). 
When this measure is unique the flow is called {\em uniquely ergodic}.
Suppose that $(X,T)$ is a uniquely ergodic flow and that the unique invariant 
measure, say $\mu$, has full support. It then follows that $(X,T)$ is minimal,
as otherwise there would be a nonempty, closed and invariant subset 
$Y \subsetneq X$ which by the 
Krylov-Bogolubov theorem will support
an invariant probability measure $\nu$ which is distinct from $\mu$ (having a smaller support).
A flow is {\em strictly ergodic} when it is uniquely ergodic and minimal.

\br

A classical theorem of Harald Bohr asserts that every {\em Bohr almost periodic function} on $\Z$
(i.e. a function in $\ell^{\infty}(\Z)$ whose orbit under translations is
norm pre-compact) can be uniformly
approximated by linear combinations of functions of the form
$$
n \mapsto e^{2\pi i n \al} \quad (n \in \Z, \al \in \R).
$$
Another characterization of almost periodic functions is the following one.
The function $f \in \ell^\infty(\Z)$ is almost periodic iff there exists 
a minimal equicontinuous flow $(X,T)$, 
a continuous function $F \in C(X)$ and a point $x_0 \in X$ with
$$
f(n) = F(T^n x_0) \quad (n \in \Z).
$$

Let $\Dcal$ be the family of functions $f \in \ell^\infty(\Z)$ such that :
there exists a minimal distal flow $(X,T)$, 
a continuous function $F \in C(X)$ and a point $x_0 \in X$ with
$$
f(n) = F(T^n x_0) \quad (n \in \Z).
$$
(We say that $f$ is {\em coming} from the flow $(X,T)$.)
Then $\Dcal$ is a uniformly closed, translation invariant subalgebra of $\ell^\infty(\Z)$,
called {\em the algebra of distal functions}.
For a fixed irrational $\al \in \R$ the function
$$
n \mapsto e^{2\pi i n^2 \al} 
$$
is distal. To see this define the flow $(X,T)$ where $X = \T^2 = (\R/\Z)^2$ (the two torus) and $T : X \to X$ 
is given by
$$
T(x,y) = (x + \al, y +2x + \al), \quad \pmod 1.
$$
Let $F(x,y)=e^{2\pi i y}$ and observe that $T^n(0,0)= (n\al, n^2 \al)$.
Thus
$$
f(n) = F(T^n(0,0)) = e^{2\pi i n^2 \al}. 
$$
It is easy to see that $(X,T)$ is a distal flow and it follows that $f$ is a distal function.

It is shown in \cite{Fur-81} that for every $\al \in \R$ and $k \in \N$ one can build a 
distal flow $(X,T)$ on some torus $X$ so that, with the right choice of $F \in C(X)$ and $x_0 \in X$, the function
$$
f_{\al,k} : n \mapsto e^{2\pi i n^k \al} \quad (n \in \Z)
$$
has the form $f_{\al,k}(n) = F(T^n x_0)$. Thus the family $\{f_{\al,k} : \al \in \R, k \in \N\}$ consists of
distal functions.

Let us call the smallest closed, translation invariant subalgebra of $\ell^\infty(\Z)$ containing
$\{f_{\al,k} : \al \in \R, k \in \N\}$, the {\em Weyl algebra} and denote it by $\Wcal$.
We have $\Wcal \subset \Dcal$. Recalling Bohr's theorem
it is natural to ask whether $\Wcal = \Dcal$ ? (To the best of my knowledge this was first
asked by P. Milnes.)
The answer to this question is negative and we will see this in two different approaches
In section \ref{sec-Wse} we will show that every function of $\Wcal$ is {\em strictly ergodic}  (i.e. 
it comes from a strictly ergodic flow).
Since in \cite{Fur-61} Furstenberg produces a minimal flow which is not strictly ergodic,
this implies that $\Wcal \not= \Dcal$. Moreover, we will show that the family 
$\Scal^d$ of strictly ergodic functions in $\Dcal$ does not form an algebra and hence in
particular does not coincide with $\Wcal$.

In section \ref{sec-se} we shall see that a {\em multiplier for strict ergodicity},
either within $\Dcal$ or in general (see definition below), is necessarily a constant. 
In Section \ref{sec-new} we give 
an example of a strictly ergodic, metric, distal flow whose 
which admits a non-strictly ergodic $2$-fold minimal self-joining.
It then follows that the enveloping group of this flow is not strictly ergodic
(as a $T$-flow).
In section \ref{sec-nil}
we show that the functions coming from a certain translation on a 3-dimensional nil manifold are
not in $\Wcal$, demonstrating again the inequality of $\Wcal$ with $\Dcal$.
In particular a classical $\Theta$-function defined on $\R^3$ yields such a concrete function.

\begin{rmk}
\begin{enumerate}
\item
Independently of us, and at about the same time, E. Salehi has also shown that
$\Wcal \subsetneq \Dcal$, \cite{S}. See also \cite{Kn-67} and \cite{JN-10}.
\item
In his paper 2016 paper \cite{R} Juho Rautio presents an in depth research of the flow $|\Wcal|$.
The main results of this definitive paper are as follows.
\begin{enumerate}
\item[(i)] 
The flow $|\Wcal|$ has quasi-discrete spectrum (see \cite{H-P}), 
and is in fact the universal quasi-discrete spectrum flow;
it admits every minimal quasi-discrete spectrum flow as a factor.
Thus the Weyl algebra $\Wcal$ coincides with the algebra generated by the minimal
systems having quasi-discrete spectrum.
\item[(ii)]
An example of a factor of the flow $|\Wcal|$ which does not have quasi-discrete spectrum is constructed.
\item[(iii)]
An explicit description of the topological and algebraic structures of the right topological group $|\Wcal|$
is given.
\end{enumerate}
\end{enumerate}
\end{rmk}

We conclude this introduction with a brief description of the interplay
between functions in $\ell^\infty(\Z)$, $T$-subalgebras of $\ell^\infty(\Z)$ and pointed flows.
(A {\em pointed flow} $(X, x_0, T)$ is a flow with a distinguished point $x_0 \in X$
whose orbit is dense, i.e. $\overline{\{T^n x_0 : n \in \Z\}} = X$.)

Starting with a function $f \in \ell^\infty(\Z)$ one can form the flow $X_f$,
the closure of the orbit of $f$ under translation in the {\em Bebutov flow} $(\R^\Z,T)$
(with product topology). We usually consider this as a pointed flow $(X, f, T)$.
Given a pointed flow $(X, x_0, T)$ we let
$$
al(X,x_0) =
\{ f \in \ell^\infty(\Z) : \exists \ F \in C(X), \ f(n) = F(T^n x_0)\ (n \in \Z)\}.
$$
This is a {\em $T$-subalgebra} of $\ell^\infty(\Z)$, i.e. a translation invariant, uniformly closed,
closed under conjugation subalgebra, which is isometrically isomorphic to $C(X)$.
The Gelfand space $|\Acal|$ of a $T$-subalgebra $\Acal$ of $\ell^\infty(\Z)$ is a flow under $T$,
the homeomorphism induced on $|\Acal|$ by translation on $\Acal$.
Again we consider $|\Acal|$ as a pointed flow where  the
multiplicative functional $f \mapsto f(0); \ \Acal \to \C$, is the distinguished point.
We then have $|al(X_f)| \cong X_f, \ al(|\Acal|) = \Acal$ and $al(X_f) = al(f)$,
where the latter is the smallest $T$-subalgebra containing $f$.

If $\{(X_i, x_i, T)\}_{i \in I}$ is a family of pointed flows then their sup :
$\bigvee_{i \in I}(X_i, x_i)$ is the subflow of the product flow $\prod_{i \in I} X_i$
which is the orbit closure of the point $x \in \prod_{i \in I} X_i$ defined by $x(i) = x_i$.
We consider this flow as a pointed flow with distinguished point $x$.
The smallest $T$-subalgebra containing $\Acal(X_i, x_i)$ \ $(i \in I)$,
is just $\Acal(\bigvee_{i \in I}(X_i, x_i), x)$. 

A function $f \in \ell^\infty(\Z)$ is called a {\em multiplier for strict ergodicity}
(within $\Dcal$) if for every strictly ergodic pointed (distal) flow $(X, x_0, T)$,
the flow $(X_f, f) \vee (X, x_0)$ is strictly ergodic (and distal).
In section \ref{sec-nil} we assume some familiarity with R. Ellis' algebraic theory of topological dynamics 
(see \cite{Ellis} or \cite{Gl1}).

I thank Professor Benjamin Weiss for reading parts of the paper and correcting
several mistakes.
 
\br

\section{$|\Wcal|$ is strictly ergodic}\label{sec-Wse}

Let $(X,T)  \overset{\pi}{\to} (Y,T)$ be a homomorphism of minimal flows. 
$\pi$ is an {\em almost periodic extension}
iff there exists a minimal flow $(\tilde{X}, T)$ and 
homomorphisms $\sig$ and $\rho$ such that the diagram
\begin{equation*}
\xymatrix
{
\tilde{X} \ar[dr]_{\rho} \ar[r]^{\sig}  &  X \ar[d]^{\pi} \\
  & Y 
}
\end{equation*}
is commutative and $\rho$ is a group extension
(i.e., there exists a compact group $K$ of automorphisms of 
$(\tilde{X},T)$ such that $(\tilde{X}/K, T) \cong (Y,T)$).

If $(X,T)  \overset{\pi}{\to} (Y,T)$ is an almost periodic extension and $\nu$ is a $T$-invariant probability 
measure on $Y$, there exists a canonically defined $T$-invariant measure $\mu$ on $X$ as follows. Let 
$\la$ be the Haar measure on $K$.  Define $\tilde{\mu}$ on $\tilde{X}$ by
$$
\int f(x)\, d\tilde{\mu}(x) = \int \int f(xk) \, d\la(k) \, d\nu(\rho(x))
$$
and put $\mu = \tilde{\mu} \circ \sig^{-1}$. (We say that $\mu$ is the {\em Haar lift} of $\nu$.)

By the distal structure theorem of Furstenberg \cite{Fur-63}, (which asserts that every
minimal distal flow can be represented as a tower of almost periodic extensions) we see, 
using the procedure described above, that every minimal distal flow carries a canonically defined 
invariant measure.

We recall the following theorem of Furstenberg \cite{Fur-61} (see also \cite[Theorem 3.30]{Gl2}):

\begin{thm}\label{unique}
Let $(X,T)  \overset{\pi}{\to} (Y,T)$ be an almost periodic extension. Let $\nu$ be $T$-invariant probability 
measure on $Y$ and let $\mu$ be its Haar lift on $X$. 
Suppose $\mu$ is ergodic, then
for every invariant probability measure $\mu_1$ on $X$ with
$\mu_1 \circ \pi^{-1} = \nu$, we have $\mu_1 = \mu$.
In particular, if $(Y,T)$ is strictly ergodic and $\mu$ is ergodic, then $(X,T)$ is strictly ergodic.
\end{thm}

\begin{cor}\label{cor-fur}
Let $(V,T)$ be a minimal distal flow, $\mu$ its canonically defined measure.
If $\mu$ is ergodic then $(X,T)$ is strictly ergodic.
\end{cor}

\begin{proof}
Transfinite induction on the Furstenberg tower of $(X,T)$ (see \cite{Fur-63}).
\end{proof}

\begin{lem}[Berg]\label{lem-berg}
Let $(X,T,\mu)$ and $(Y,T,\nu)$ be measure preserving ergodic processes.
Suppose $(\tilde{X}, T, \tilde{\mu})$ and $(\tilde{Y}, T, \tilde{\nu})$ are disjoint,
where the latter are the largest Kronecker factors of $(X,T,\mu)$ and $(Y,T,\nu)$
respectively. Then $(X \times Y, T \times T, \mu \times \nu)$ is ergodic.
\end{lem}

\begin{cor}\label{cor-berg}
Let $(X,T,\mu)$ and $(Y,T,\nu)$ be strictly ergodic distal flows. 
Suppose that as measure preserving processes the corresponding largest
Kronecker factors $(\tilde{X}, T, \tilde{\mu})$ and $(\tilde{Y}, T, \tilde{\nu})$ are disjoint,
then $(X \times Y, T \times T, \mu \times \nu)$ is strictly ergodic.
\end{cor}

\begin{proof}
The disjointness of the measure preserving processes $(\tilde{X}, T, \tilde{\mu})$ and 
$(\tilde{Y}, T, \tilde{\nu})$ implies that the largest equicontinuous factors  of $(X,T)$ and $(Y,T)$ are
disjoint. This implies that $(X,T)$ and $(Y,T)$ are topologically disjoint and it follows that
$(X \times Y, T \times T, \mu \times \nu)$ is a minimal distal flow.
Obviously $\mu \times \nu$ is the canonical measure on this distal flow and the corollary
follows now from Lemma \ref{lem-berg} and Corollary \ref{cor-fur}. 
\end{proof}

\begin{lem}\label{lem-se}
For a fixed $\beta \in \R$ and $m \in \N$ consider the function
$$
f(n) = e^{2\pi i \beta (n + n^2 + \cdots + n^m)} 
$$
and let $X_f$ be the associated flow. Then
\begin{enumerate}
\item[(i)] $X_f$ is strictly ergodic.
\item[(ii)] The functions $e^{2\pi i \beta   n^j} \  (j=1,2,\dots,m)$ are all in $al(X_f)$.
\item[(iii)] If $\beta$ is rational $X_f$ is finite hence Kronecker.
Otherwise $(\T, ,T_\beta)$ is the largest Kronecker factor of $(X_f, T, \mu)$, where $\mu$ is the unique
$T$-invariant measure on $X_f$.
\end{enumerate}
\end{lem}

\begin{proof}
We prove the case $m = 3$. (This will indicate the proof in the general case which is cumbersome to write.)

In this case
\begin{gather*}
f(n) = \exp(2\pi i \beta (n + n^2 +  n^3)) \\
f(n+k) = \exp(2\pi i \beta ((n  +k)  + (n +k)^2 +  (n +k)^3).
\end{gather*}
Let $k_j \beta \to x, \ k_j^2 \beta \to y, \ k_j^3\beta \to z$ for some sequence $k_j \to \infty$. Then
\begin{align*}
\lim_{j \to \infty} f(n + k_j) & = g_{x, y, z}(n) \\
& =  \exp(2\pi i (n(3y + 2x + \beta) + n^2(3x + \beta) + n^3 \beta + x + y +z)),
\end{align*}
$g_{x, y,z} \in X_f$ and
\begin{align*}
T g_{x,y,z}(n) & = g_{x,y,z}(n+1)\\
& =  \exp\{2\pi i (n[3(y + 2x + \beta) + 2(x + \beta) + \beta] +
n^2[(3x + \beta) + \beta] + n^3 \beta  \\
& \ \ \ + (x + \beta) + ( y + 2x + \beta)  + (z + 3x + 3y + \beta))\}\\
& = g_{x+\beta,  y + 2x + \beta, z + 3x + 3y + \beta}(n).
\end{align*}
Thus $(X_f, T)$ can be identified with the orbit closure of $(0, 0, 0) \in \T^3$ 
under the transformation
$$
T(x,y,z) = (x +\beta, y +2x  + \beta, z + 3x +3y + \beta). 
$$
When $\beta$ is irrational $X_f = \T^3$ and $(X_f,T)$ is strictly ergodic by \cite{Fur-61}.

Now $T^n(0,0,0)= (n \beta, n^2 \beta, n^3 \beta)$ and if we let $\pi_j : \T^3 \to \T$ be the projection on the $j$th coordinate $(j = 1.2.3)$, we have
$$
\exp(2 \pi i(\pi_j \circ T^n(0,0,0))) = \exp(2\pi in^j \beta) \in al(X_f).
$$ 
By comparing Fourier expansions one can show that $f(T(x,y,z)) = \la f(x,y,z)$
with $|\la| =1, \ f \in C(X_f)$ implies $f(x,y,z) = f(x)$ and $\la = e^{2 \pi i \theta}$, where
$\theta \in \Z \beta$. This proves (iii).
\end{proof}

\begin{thm}\label{Wse}
The flow $(|\Wcal|, T)$ is strictly ergodic. Hence every function in $\Wcal$
is strictly ergodic; 
i.e., $f \in \Wcal$ implies $(X_f,T)$ is strictly ergodic.
\end{thm}
\begin{proof}
Let $\Wcal_0$ be the algebra generated by the functions
$$
\{\exp(2 \pi i n^k \al) : \al \in \R, \ k \in \Z\}.
$$
An element $f \in \Wcal_0$ has the form
$$
f(n) = \sum_{j=1}^m a_j e^{2\pi i p_j(n)},
$$
where $p_j(n) = \sum_{l =0}^{v_j}
\al_{l, j} n^l$.
We construct a strictly ergodic flow $(X,T)$ such that for some $x_0 \in X$ and for each $j$ and $l$,
$\exp(2\pi \al_{j,l}n^l)$ is in $al(X,x_0)$.

Choose a finite $\Q$-independet set $\{\beta_1, \beta_2, \dots, \beta_u\} \subset [0,1]$ such that
$\{\al_{j,l}\}_{j,l}^{m,v_j}$ is contained in $\oplus_{s =1}^u \Z \beta_s$.

Let 
$$
f_s(n) = \exp(2\pi i \beta_s(n + n^2 + \cdots + n^{L_s})) \quad (s =1,2,\ldots, u).
$$
By Lemma \ref{lem-se} the functions
$\exp(2\pi i \beta_s n^l) \ (1 \leq l \leq L)$ are
in $al(X_{f_s})$, $(X_{f_s},T)$ is strictly ergodic and the $\beta_s$ rotation is its largest Kronecker factor.
By Corollary \ref{cor-berg} $X = \prod_{s =1}^u X_{f_s}$ is also strictly ergodic and 
$f \in al(X)$. 

Thus $f \in \Wcal_0$ implies $f$ is strictly ergodic. Since 
$$
|\Wcal| = \lim_{\rightarrow} \{ |X_f| : f \in \Wcal_0\},
$$
it follows that $|\Wcal|$ is strictly ergodic and our theorem follows.
\end{proof}

Now, in view of Theorem \ref{Wse},
and the example of Furstenberg in \cite{Fur-61} of a non-uniquely ergodic minimal distal flow on $\T^2$, we obtain  the following corollary.

\begin{cor}
The inclusion $\mathcal{W} \subset \mathcal{D}$ is proper. More specifically,
there exists a metric minimal distal flow $(X,T)$ which is not strictly ergodic and hence is not in $\mathcal{W}$. 
\end{cor}

{\bf The following example shows that the family $\Scal^d$ of distal, strictly ergodic functions
does not form an algebra (or even a vector space).
In particular $\Scal^d$ contains $\Wcal$ but does not coincide with it.}

As in \cite{Fur-61} let $v_1 =1, \ v_{k+1} = 2^{v_k} + v_k + 1, \
n_k = 2^{v_k}, \ n_{-k} =- n_k $ and $ \al = \sum_{k=1}^\infty n_k^{-1}$. It was shown in 
\cite{C-G} that for some $t, \ 0 \leq t \leq 1$ the function
$$
h(x) = t\sum_{k \not= 0}(e^{2\pi i n_k \al} -1) e^{2\pi i n_k x}
$$
is distal and strictly ergodic and by the same argument so is the function
$$
g(x) = t\sum_{k \not= 0}(e^{2\pi i n_k \al} -1)(1 + \frac{1}{|k|})e^{2\pi i n_k x}.
$$
However the function 
$$
f(x) =  g(x) - h(x) = t\sum_{k \not= 0} \frac{1}{|k|}(e^{2\pi i n_k \al} - 1)e^{2\pi i n_k x}
$$
is not strictly ergodic by \cite{Fur-61}.

\section{Multipliers for strict ergodicity are constants}\label{sec-se}

In order to motivate our next problem we survey quickly the following well known results from topological dynamics.

Let us call a function $f \in \ell^\infty(\Z)$ {\em minimal} if there is
a minimal flow $(X,T)$, a point $x_0 \in X$ and a function $F \in C(X)$ with
$$
f(n) = F(T^n x_0).
$$
Let $\Mcal$ be the family of of all minimal functions in $\ell^\infty(\Z)$. It is well known that
$\Mcal$ is not preserved under sums. Using Zoren's lemma  we can however form the family $\{\Mcal_\al\}$ of maximal $T$-subalgebras of $\Mcal$. These are pairwise isometric (as Banach algebras) and in fact there exists
 a universal minimal flow $(M,T)$  such that each $\Mcal_\al$ is of the form $al(M,m)$ for some $m \in M$.
 The object we are interested in now is the $T$-algebra $\Lcal = \cap \Mcal_\al$. One can easily see
 that $f \in \Lcal$ iff the corresponding
 minimal pointed flow $(X_f,f)$ has the following propery: For every minimal pointed flow $(Y,y_0)$
 the flow $(X_f \vee Y, (f,y_0))$ is minimal. We say that $f$ is a {\em multiplier for minimality}.
 
 One more definition. A minimal flow $(X,T)$ is called {\em point distal}
 if there exists a point $x_0 \in X$ such that no point of $X$ other than $x_0$ is proximal to $x_0$.
 We call $f \in \ell^\infty(\Z)$ {\em point distal}
 if it is coming from a point distal flow; i.e. if there is a point distal flow $(X,T)$ and a point $x_0 \in X$
 as above with 
 $$
 f(n) = F(T^n x_0)
 $$
for some $F \in C(X)$. We now have the following surprising 
characterization of multipliers for minimality, \cite[Theorem III.8]{AH}.

\begin{thm}\label{thm-pointdistal}
$\Lcal = \cap \Mcal_\al$ coincides with the algebra of point distal functions.
\end{thm}

For later use we also need to mention the following.

\begin{thm}\label{disj}
Every minimal weakly mixing flow is disjoint from every point distal flow.
\end{thm}

(A flow $(X,T)$ is {\em weakly mixing} if $(X \times X, T \times T)$ is topologically ergodic,
i.e., an open invariant non-empty sunset of $X \times X$ is necessarily dense. It is {\em strongly mixing} if
for every non-empty open sets $A$ and $B$ the set $N(A,B) = \{n \in \Z : T^n A \cap B \not= \emptyset\}$ has a finite complement.)

Now back to strictly ergodic functions. We have seen that the set $\Scal^d$ of strictly ergodic distal functions
does not form an algebra. The same example shows that neither is the family $\Scal$ 
of strictly ergodic functions an algebra. Using Zoren' lemma again we can form the two families
$\{\Scal_\beta\}$ and $\{\Scal_\ga^d\}$ of maximal $T$-subalgebras in $\Scal$ and $\Scal^d$, respectively.
Motivated by the example of point distal functions we pose the problem of identifying 
$\cap \Scal_\beta = \Scal_0$ and $\cap \Scal_\gamma^d = \Scal_0^d$.

As in the case of $\cap \Mcal_\al$ it is easy to see that being an element of $\cap \Scal_\beta $
is the same as being a multiplier for strict ergodicity. I.e. $f \in  \cap \Scal_\beta$ iff 
for every strictly ergodic flow $(Y,T)$ and every point $y \in Y$ the flow $(X_f \vee Y, (f,y))$ is strictly ergodic.
A similar statement holds for $\cap \Scal_\gamma^d$.

The answer to our question is disappointing.

\begin{thm}\label{trivial}
The only multipliers for both strict ergodicity and 
strict ergodicity within $\Dcal$, are the constant functions: i.e. $\Scal_0 = \Scal_0^d = \C$.
\end{thm}

Notice that $\Scal_0 = \cap \Scal_\beta = \C$ does not automatically imply
$\Scal^d_0 = \cap \Scal_\ga^d = \C$; it is not at all clear that for a given $\beta$
 the algebra $\Scal_\beta \cap \Dcal$ is not a proper subalgebra of one of the maximal
algebras $\Scal_\ga^d$.

For the proof we need the following strengthening of the Jewett Krieger theorem due to Lehrer \cite{Leh}.

\begin{thm}\label{lehr}
If $(\Om, \Fcal,\mu,T)$ is a properly ergodic process then there exists a strictly ergodic
flow $(X,T)$ with unique invariant measure $\nu$ such that $(Y,T)$ is  topologically strongly mixing
and $(Y, \Bcal, \nu, T)$ ($\Bcal =$ Borel field of $Y$) is measurably isomorphic to $(\Om, \Fcal, \mu, T)$.
\end{thm}

For the distal case we need the following theorem which can be deduced from \cite{Fur-61}.

For $\beta \in \R$ we let $(X_\beta, R_\beta)$ be the rotation by 
$\beta$ of $\T = X_\beta$ when $\beta$ is irrational, and
$\{1,2,\ldots,n\} = X_\beta$ when $\beta$ ia a rational number of order $n$ (considered as an element of $\T$).

\begin{thm}\label{ab}
Let $\beta \in \R$ be given, then there exist an irrational $\al \in \R$ such that $\al$ and $\beta$
are rationally independent, and a strictly ergodic distal flow $(X,T)$ with invariant measure $\mu$ such that
\begin{enumerate}
\item[(i)] 
$(X, \mu,T)$ is measure theoretically isomorphic with $(X_\al \times X_\beta, R_\al \times R_\beta)$.
\item[(ii)]
There exists a continuous homomorphism $(X,T) \to (X_\al,R_\al)$.
\item[(iii)]
The flows $(X,T)$ and $(\T, R_\beta)$ are topologically disjoint.
\end{enumerate}
\end{thm}

\begin{proof}
We go back to the construction described by Furstenberg \cite{Fur-61} on page 385.
As in \cite{Fur-61} let $v_1 =1, \ v_{k+1} = 2^{v_k} + v_k + 1, \
n_k = 2^{v_k}, \ n_{-k} =- n_k $ and $ \al = \sum_{k=1}^\infty n_k^{-1}$.
We have
\begin{equation}\label{ineq}
|n_k \al - [n_k \al]| < \frac{2 \cdot 2^{v_k}}{2^{v_{k+1}}} = 2^{-n_k}.
\end{equation}

Now set
$$
h(x) = \sum_{k \not= 0}\frac{1}{|k|}(e^{2\pi i n_k \al} -1) e^{2\pi i n_k x}
$$
and let 
$$
g(e^{2\pi i x}) = e^{2\pi i t h(x)},
$$
where $t$ is yet undetermined.
By (\ref{ineq}) $h(x)$ and therefore $g(\zeta)$ are $C^{\infty}$ functions.
Now we have $h(x) = H(x+\al) - H(x)$ where
$$
H(x) = \sum_{k \not= 0}\frac{1}{|k|}e^{2\pi i n_kx}
$$
As in \cite{Fur-61} we observe that $H(x)$ is in $L^2(0,1)$, hence defines a measurable function,
which however, can not correspond to a continuous function,
hence not summable at $\theta =0$ (see  e.g. \cite[Theorem 3.4. page 89]{Z}). 

We conclude that there is some $t$ such that, {\bf for all $k \in \N$} the function
$e^{2 \pi i k t H(x)}$ will not be a continuous function (see \cite[Proposition A1, page 83]{EH}).
Taking $R(e^{2 \pi i x}) = e^{2 \pi i t H(x)}$ we have 
\begin{equation*}\label{chom}
R(e^{2 \pi i \al} \zeta) / R(\zeta) = g(\zeta)
\end{equation*}
with $R(\zeta)$ measurable but not continuous.
Now consider the transformation on $X = \T^2$ 
$$
T(x, y) = T_\phi(x,y) = (x + \al, y +t h(x) + \beta), 
$$
defined by means of the cocycle $\phi(x) : = t h(x) + \beta$
\footnote{If the $\al$ we just constructed and the given $\beta$ happen to be  rationally dependent, we can slightly modify the construction of $\al$ so that it becomes independent of $\beta$.}.

The claim (i) of the theorem follows because
the function 
$\phi(x) = th(x) + \beta = tH(x +\al) - tH(x) + \beta$ is co-homologous to $\beta$.

Claim (ii) is clear.

Thus the measure theoretical system $(X, \mu, T)$ is ergodic and by a theorem of Furstenberg
(see e.g. \cite[Theorem 4.3]{Gl2})
the flow $(X,T)$, which is a $\T$-extension of the strictly ergodic flow $(\T,R_\al)$, is also strictly ergodic.

Next observe that there is no nontrivial intermediate extension
$$
(X,T) \to (Y,T) \overset{\rho}{\to} (\T, R_\al),
$$
with $\rho$ finite to one.
In fact, any intermediate flow $(Y,T)$ is obtained as a quotient $Y \cong X / K$,
where $K$ is a closed subgroup of the group $\mathfrak{A}: = \{A_d : d \in \T\}$ of automorphisms
of the flow $(X,T)$ of the form $A_d(x,y) = (x, y+d)$.
But, any proper closed subgroup $K$ of $\mathfrak{A}$ is finite, so that
$Y \cong X/K  \overset{\rho}{\to}  (\T, R_\al)$ can not be finite to one.

%

Now in this situation one can check that the rotation $(\T, R_\al)$ is the maximal equicontinuous factor of
$(X,T)$ iff 
\begin{quote}
For every $\la \in \C$, for every $0 \neq k \in \Z$ the functional equation
$$
f(x + \al) e^{2\pi i [k(t h(x) + \beta)]}= \la f(x)
$$ 
has non non-zero continuous solution.
\end{quote}

Suppose $f$ is such a non-zero solution. Then
$$
\frac{f(x + \al)}{f(x)} \cdot \frac{e^{2\pi i k t H(x + \al)}}{e^{2\pi i kt H(x)}} = \la  e^{- 2\pi i k \beta}.
$$
Writing $F(x) =  e^{2\pi i kt H(x)}$ and $b =e^{2\pi i k \beta}$,   we have
$$
\frac{f(x + \al)}{f(x)} \cdot \frac{F(x + \al)}{F(x)} = \la  b^{-1}.
$$
Thus the function $f\cdot F$ is an eigenfunction of the flow $(\T, R_\al)$,
hence, for some $l \in \Z$, $\la b^{-1} = l \al$, and  $f \cdot F$ has the form $e^{2\pi i l x}$.
This  contradicts the fact that $F$ is not equal a.e. to a continuous function.

Therefore $(\T, R_\al)$ is the maximal equicontinuous factor of $(X,T)$ and
it follows that indeed $(X,T)$ and $(\T, R_\beta)$ are topologically disjoint.
This completes the proof of part (iii). 
\end{proof}

We are now ready for the proof of Theorem \ref{trivial}.

\begin{proof}[Proof of Theorem \ref{trivial}]

\ 

(i) $\Scal_0 = \C$.

Since every strictly ergodic flow is minimal we have for each $\beta$ an $\al$ such that
$\Scal_\beta \subset M_\al$. On the other hand, given n $\Mcal_\al$ 
and a strictly ergodic flow $(X,T)$, there is a point $x \in X$ such that 
$al(X,x) \subset \Mcal_\al$, so that each $\Mcal_\al$ contains exactly one $\Scal_\beta$.
It follows that
$$
\Scal_0 = \cap \Scal_\beta \subset \cap \Mcal_\al = \Lcal.
$$
\end{proof}

Now fix an element $ f \in \Scal_0$, and let $(X,x_0) = (X_f, f)$. Since $f \in \Lcal$
it follows that $(X, x_0)$ is point distal. Since $f \in \Scal_0$ the flow $(X,,T)$  is strictly ergodic with 
a unique invariant measure $\mu$. 
Use Lehrer's theorem, Theorem  \ref{lehr}, to produce a strictly ergodic, topologically strongly mixing flow $(Y,T)$, with invariant measure $\nu$, 
so that $(Y, \nu, T)$ is measure theoretically isomorphic to $(X, \mu, T)$; say $X \overset{\phi}{\to} Y$.

By Theorem \ref{disj} the flow $(X \times Y, T \times T)$ is minimal hence
$(X,x_0)\vee (Y, y_0) = X \times Y$
for every choice of $y_0 \in Y$. Now on $X \times Y$ 
we have the following $T \times T$-invariant measures
$$
\mu \times \nu = \int (\del_x \times \nu)\, d\mu(x)
$$
and
$$
 \int (\del_x \times \del_{\phi(x)})\, d\mu(x).
$$
By the uniqueness of the disintegration of measures these two invariant measures coincide iff
$\del_{\phi(x)} = \nu$ a.e., iff $\nu $ is a point mass, iff $(X,T)$ is a trivial one point flow,
iff $f$ is a constant.
But by assumption $f$ is a multiplier for strict ergodicity so that
$(X,x_0,T)\vee (Y, y_0,T) = (X \times Y, T \times T)$ is strictly ergodic and we conclude that $f$ is a constant.

(ii) $\Scal_0^d = \C$.

Let $\beta \in \R$ be given and let $\al$ and $(X,T)$ be as in Theorem \ref{ab}. 
Since $(X_\al, R_\al)$ and $(X_\beta, R_\beta)$ are topologically disjoint so are 
$(X,T)$ and $(X_\beta,R_\beta)$. 
This follows since by Theorem \ref{ab} $(X_\beta, \R_\beta)$ is the maximal 
Kronecker factor of $(X,T)$, and e.g., by  Theorem 4.2 in \cite{EGS}.
Thus
$(X \times X_\beta, T \times R_\beta) = (X,x_0,T)\vee (X_\beta, x_\beta, R_\beta)$ is minimal.

Now for $f(n) = e^{2\pi i n \beta}$ we have $|al(f)| = (X_\beta, R_\beta)$ and if $f \in \Scal_0^d$
then $f$ is a multiplier for strict ergodicity within $\Dcal$ and $(X \times X_\beta, T \times R_\beta)$
is strictly ergodic. However, if $\mu$ 
and $\nu$ are the invariant measures on $X$ and $X_\beta$ respectively and $X \overset{\phi}{\to} X_\beta$
is the measure theoretical isomorphism, we can disintegrate $\mu$ over $\nu$
$$
\mu = \int_{X_\beta} \mu_{z} \, d\nu(z)
$$
($\mu_z$ a probability measure on $\phi^{-1}(z)$ for a.e. $z \in X_\beta$)
and then $\mu \times \nu$ and 
$\int (\mu_z \times \nu) \, d\nu(z)$ are two
different invariant measures  on $(X \times X_\beta, T \times R_\beta)$, a contradiction.

Thus for every $\beta \in \R$, $f_\beta(n) = e^{2\pi i n \beta}$ 
is not in $\Scal_0^d$ and since $\Scal_0^d$ is a distal algebra 
it follows from the structure theorem  for minimal distal flows \cite{Fur-63} that $\Scal_0^d = \C$.

\begin{rmk}
It is interesting to compare the discussion above, about multipliers for strict ergodicity,
with the study of multipliers for minimal weakly mixing flows in \cite{EG}.
There it was shown that the collection of minimal functions which
are multipliers for weak mixing, which by definition is the
intersection of the family $\{\mathcal{A}_\del\}$ of maximal
weakly mixing subalgebras of $\ell^{\infty}(\Z)$ (sitting in the algebra which corresponds to the
algebra $C(M,u)$, where $M$ is the universal minimal flow), 
coincides with the algebra of {\em purely weakly mixing functions} hence, in particular, is nontrivial.
\end{rmk}

\section{A strictly ergodic distal flow with a non strictly ergodic self-joining}\label{sec-new}

It is not hard to see that the family of Weyl minimal flows is closed under self joinings.
As a consequence the enveloping group of a Weyl flow is again a Weyl flow, hence strictly ergodic.
The same assertions hold for the family of pro-nil flows
(i.e. those flows that correspond to subalgebras of the $T$-algebra generated by all the nil-functions;
see Section \ref{sec-nil} below).  
Is it true that the enveloping group of every strictly ergodic distal flow is strictly ergodic ?

Our goal in this section is to construct a strictly ergodic distal flow with a non strictly ergodic self-joining.
The enveloping group of such a flow is clearly not strictly ergodic.

We start with the Anzai-Furstenberg example, the simplest example of a minimal 
distal but not equicontinuous flow.
On the two torus $\T^2$ we let
$$
A(x,z) = (x + \al, z +x)
$$
where $\al$ is irrational. We observe that for any $\beta \in \R$ 
the self joining $(X,(0,0),A) \vee (X,(\beta,0),A)$ admits the flow $(\T, R_\beta)$ as a factor,
with factor map
$$
((x,z)(x',z')) \mapsto z' -z.
$$

Next consider the strictly ergodic distal flow $(X,T)= (\T^2, T_\phi)$ 
constructed in Section \ref{sec-se} with an irrational $\beta$ independent of $\al$, 
and define $S : \T^3 \to \T^3$ by
$$
S (x,y,z) = (x +\al, y + \phi(x), z + x).
$$

\begin{claim}
The flow 
$$
(W,S) := (\T^3,S) \cong (X, (0,0), T) \vee (\T^2,(0,0), A) = (X \underset{(\T, R_\al)}{\times}\T^2, T \times A)
$$ 
is strictly ergodic.
\end{claim}

\begin{proof}

Let $\la$ be Lebesgue measure on $\T^3$.
If $f(x,y,z)$ is an $S$-invariant function in $L^2(\la)$ then using the Fourier expansion
$$
f(x,y,z) = \sum_{k \in \Z} a_k(x,y) e^{2\pi i kz},
$$
we can check that $f$ is a constant, as follows.

Since $f$ is invariant we have, for a.e. $(x, y, z) \in \T^3$,
$$
f(x + \al,y + \phi(x),z + x) = \sum_{k \in \Z} a_k(x + \al ,y + \phi(x))e^{2\pi i kx} e^{2\pi i kz},
$$
Thus for $k \neq 0$, a.e.
$$
a_k(x,y) = a_k(x + \al ,y + \phi(x))e^{2\pi i kx}.  
$$
By ergodicity of $(X,T_\phi)$ either $a_k = 0$ a.e. or
it vanishes on a set of measure zero. 
Therefore we get an a.e. equality
$$
\frac{a_k(T(x,y))}{a_k(x,y)} = e^{-2\pi i kx}.
$$

However, via the the map $J : \T^2 \to \T^2, \ J(x,y) = (x, y + tH(x))$ the flow 
$(X,T_\phi)$ is measure theoretically isomorphic to the 
flow 
$$
T_{\al, \beta} : \T^2 \to \T^2, \quad  T_{\al,\beta}(x,y) = (x + \al, y +\beta)
$$
and thus, with $A_k(x,y) : = a_k \circ J (x,y)$, we get 
\begin{align*}
A_k(x,y) & = (a_k \circ J)(x,y) = a_k(x, y +tH(x))\\
& = a_k(x +\al, y + tH(x) + th(x) + \beta) e^{2\pi i kx}\\
&= a_k(x +\al, y +t H(x+\al) + \beta) e^{2\pi i kx}\\
&= (a_k \circ J)(x +\al, y +\beta)e^{2\pi i kx} = A_k(x +\al, y + \beta)e^{2\pi i kx}
\end{align*}\
hence
$$
\frac{A_k(x +\al,y + \beta)}{A_k(x,y)} = e^{-2\pi i kx}.
$$
Using again Fourier decomposition for $A_k$ we conclude 
that for $k \neq 0$ 
we have $A_k(x,y) =0$, hence also $a_k =0$. 
Indeed, with
$$
A_k(x,y) = \sum_{l \in \Z} b_l(x) e^{2\pi i lky},
$$
we get
$$
A_k(x,y)  = \sum_{l \in \Z} b_l(x) e^{2\pi i lky} =
\sum_{l \in \Z} b_l(x +\al) e^{2\pi i lky} e^{2\pi i lk \beta} e^{2\pi i kx},
$$
hence, for $l \neq 0$, 
$$
b_l(x) =b_l(x +\al) e^{2\pi i lk \beta} e^{2\pi i kx}.
$$

Next let $b_l(x) = \sum_{m \in \Z} d_n e^{2\pi i m x}$. 
Then, comparing coefficients we get,
$d_m = d_{m - k}e^{2 \pi i (lk\beta+ (m -k)\al)}$, 
hence $d_m = 0$, hence $b_l(x) =0$.
Thus $A_k(x,y) =  b_0(x)$, hence $b_0(x) = b_0(x +\al) e^{2 \pi i kx}$
and again Fourier expansion shows that $b_0 =0$, and consequently $A_k =0$.

Finally for $k =0$ the equation we get is
$$
\frac{A_0(x +\al,y + \beta)}{A_0(x,y)} = 1,
$$
hence $A_0(x,y)$ is a constant and so is $f$.

Thus the measure theoretical system $(W, \la, S)$ is ergodic and by a theorem of Furstenberg
(see Theorem \ref{unique})
the flow $(W,S)$, which is a $\T$-extension of the strictly ergodic flow $(X,T)$, is also strictly ergodic.
\end{proof}

We will now show that the self-joining 
$$
M : = (W,(0,0,0),S) \vee (W, (\beta, 0,0),S)
$$ 
is not strictly ergodic.
In fact, since both $(X,T)$ and $(\T, R_\beta)$ are  factors of $M$, it follows that also
the flow $(X,x_0,T)\vee (\T, 0, R_\beta)$
(with $x_0 = (0,0,0)$)  is a factor of $M$.
However, as we have seen in the proof of Theorem \ref{trivial}, this flow
$ (X,x_0,T)\vee (\T, 0, R_\beta) = (X \times \T, T \times R_\beta)$ 
is minimal but not strictly ergodic. Thus we conclude that also our minimal flow $M$ 
is not strictly ergodic.

\begin{cor}
There exists a metric, strictly ergodic, distal flow $(W,T)$ whose enveloping group 
$E(W,T)$ is not strictly ergodic.
\end{cor}

\begin{proof}
Recall that, as flows,
$$
E(W,T) \cong \bigvee_{w \in W} (W,w,T).
$$
It follows that in the example above, $M$ is a factor of $E(W,T)$.
\end{proof}


%

\section{A nil-flow relatively disjoint from $\Wcal$}\label{sec-nil}

Our goal in this section is to show that the 3-dimensional nil-flow
(by this we mean the time one, discrete flow in a one parameter transformation group), 
which is a strictly ergodic and distal flow, is relatively disjoint over its equicontinuous factor from $|\Wcal|$.
Thus a function $f \in \ell^\infty(\Z)$ coming from this 
nil-flow - $f$ not almost periodic - will be an element of $\Scal^d$ but not an element of $\Wcal$.
One such example is the $\Theta$-function:
$$
f(n) = F(n\al, n\beta, \frac{\al\beta}{2} + n \ga),
$$
where 
$$
F(x,y,z) = e^{2\pi i z} \sum_{m \in \Z} e^{2 \pi i m x} e^{-\pi(m+y)^2}
$$
and $\al,  \beta, \ga \in \R$, $\al, \beta$ irrationals, rationally independent. 

Let
$$
N=\biggl\lbrace\begin{pmatrix} 1&x&z\\
0&1&y\\
0&0&1\end{pmatrix} :x,y,z\in\mathbb {R}\bigg\}\ ,
$$
and $\Ga \subset N$ the discrete co-compact subgroup with integral coefficients.
Put $Z = N/\Ga$ and let
$$
T=\begin{pmatrix}  1&\al&\ga + \frac12 (\al \beta) \\
0&1&\beta\\
0&0&1\end{pmatrix} ,
$$
It was shown in \cite{AGH} that the flow $(Z,T)$ (where $T(g\Ga) = (Tg) \Ga \ (g \in N)$) is 
strictly ergodic, distal and not equicontinuous; in fact the largest equicontinuous factor of $(Z,T)$
is the flow
$(Y,T) = (\T^2, R_\al \times R_\beta) \cong (N/B\Ga, T)$ where
$$
B=\biggl\lbrace\begin{pmatrix}  1& 0 &z \\
0&1&0\\
0&0&1\end{pmatrix}
: z \in \R\bigg\}.
$$
We let $(Z,T) \overset{\pi}{\to} (Y,T)$ be the canonical homomorphism.
In the next lemma we use the notations introduced in the proof of Theorem \ref{Wse}. Thus we let
$$
X = X(\beta_1, \dots, \beta_u; L_1, \ldots, L_u) = \prod_{s =1}^u X_{f_s},
$$
where 
$$
f_s(n) = \exp(2\pi i \beta_s(n + n^2 + \cdots + n^{L_s})),
$$
$(s =1,2,\ldots, u,\  n \in \Z )$, $\{\beta_1, \ldots, \beta_u\}$ is a $\Q$-independent subset of $\R$
and $L_1, \ldots, L_u$ are arbitrary positive integers.

We have seen in the proof of Theorem \ref{Wse} that $|\Wcal|$ is the inverse limit of the flows
$X = X(\beta_1, \dots, \beta_u; L_1, \ldots, L_u)$. We observe that topologically each such $X$ has the form
$\T^\nu \times K$ where $\nu$ is a positive integer and $K$ is a finite set. 
Moreover, if $\nu >1$ then  $X$ is a $\T$-extension  of a flow $X'$ on $\T^{\nu-1} \times K$ of the same form.
Also note that $\pi_1(X) = \Z^\nu$, where $\pi_1$ denotes fundamental group.

\begin{lem}\label{lem-disj}
Let $X = X(\beta_1, \dots, \beta_u; L_1, \ldots, L_u)$, where $\beta_1 = \al, \ \beta_2 = \beta$.
Then $(Y,T)$ is also a factor of $(X,T)$ and
$(X,T)$ is relatively disjoint from $(Z,T)$ (the nil-flow with parameters $\al, \beta, \ga$)
over their common equicontinuous factor $(Y,T)$; i.e. the flow
$(X \underset{Y}{\times} Z, T \times T)$ is minimal, where
$$
X \underset{Y}{\times} Z = \{(x, z) : \phi(x) = \pi(z)\}
$$
and $\phi$ and $\pi$ are the homomorphisms of $X$ and $Z$, 
respectively, onto $(Y,T) \cong (\T^2, R_\al \times R_\beta)$.
\end{lem}

\begin{proof}
Choose $x_0 \in X, \ z_0 \in Z$, with $\phi(x_0) = \pi(z_0) = y_0$ and let
$$
(W, w_0, T) = (X, x_0, T) \vee (Z, z_0, T) \subset X \underset{Y}{\times} Z.
$$ 
$(W,T)$ is a distal minimal flow. We let $W \overset{\psi}{\to} Z$ and $W \overset{\eta}{\to} X$ be
the projections.
If $(x, z), (x, z') \in \eta^{-1}(x)$ then by the commutativity of the diagram
\begin{equation*}
\xymatrix
{
& (W,T) \ar[dl]_\eta \ar[dr]^\psi & \\
(X,T) \ar[dr]_\phi & & (Z,T) \ar[dl]^\pi\\
& (Y,T) &
}
\end{equation*}
$\pi(z) = \pi(z')$.
By considering the relative regionally proximal relation $Q_\eta$ it is now clear that $\eta$
is an almost periodic extension.
We let $F = \mathfrak{G}(Y,y_0), \ A = \mathfrak{G}(X,x_0), \ B = \mathfrak{G}(Z, z_0)$ be the
corresponding Ellis groups.
We have $B \nor F$ with $F/B \cong \T$ and $\mathfrak{G}(W,w_0) = A \cap B : = D$.
The fact that $\eta$ is an almost periodic extension is equivalent to $A/D$ being a 
compact Hausdorff space with respect to the its $\tau$-topology (see e.g. \cite{Gl1}).
Put $D_0 = \cap_{\al \in A} \al D \al^{-1}$,
then $D_0 \nor A$ and $K = A/D_0$ is a compact Hausdorff topological group.

Now our assertion that $X \underset{Y}{\times}Z$ is minimal is equivalent to $AB =F$. 
If the latter equality does not hold then the group $AB/B \subset F/B \cong \T$, 
is a proper closed subgroup of $\T$ hence finite. Consider the map 
$$
aD_0 \mapsto aB : A/D_0 \to F/B.
$$ 
Since $D_0 \subset B$ this is an isomorphism of the group $A/D_0$ onto the finite group
$AB/B$. Thus $A/D_0$ and a fortiori $A/D$ are finite sets. 
Since $(W,T)$ is distal this means that  $(W,T) \overset{\psi}{\to} (Z,T)$
is finite to one. Hence $\psi$ is a covering map and $\eta_* : \pi_1(W, w_0) \to \pi_1(Z,z_0)$
is a monomorphism.
We conclude that $\pi_1(W,w_0)$ is a subgroup of $\Z^\nu$ and in particular abelian.

Let $(X,T) \overset{\theta}{\to}(X',T)$ be the $\T$-extension discussed above
(just before Lemma \ref{lem-disj}).
This fits into the following commutative diagram
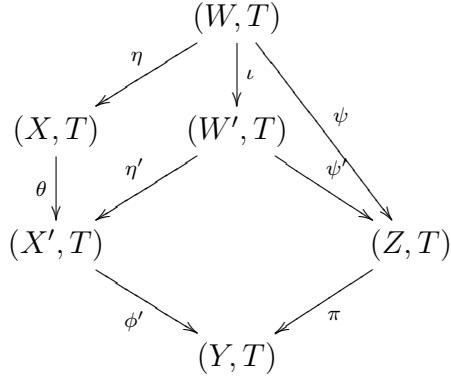
\begin{equation*}
\xymatrix
{
& (W,T)\ar[dl]_{\eta} \ar[d]^{\iota} \ar[ddr]^{\psi}& \\
(X,T) \ar[d]_\theta & (W',T) \ar[dl]_{\eta'} \ar[dr]^{\psi'} & \\
(X',T) \ar[dr]_{\phi'} & & (Z,T) \ar[dl]^\pi\\
& (Y,T) &
}
\end{equation*}
where by induction hypothesis (on $\nu$) the flow
$$
(W',T) = X' \underset{Y}{\times}Z
$$ 
is minimal.

Now $W  \overset{\iota}{\to} W'$ is a group extension. If the corresponding fiber group is $\T$, it follows that
$(X \underset{Y}{\times}Z,T)$ is minimal and we are done. Otherwise $\iota$ is a finite to one extension
and again we have that 
$\iota_* : \pi_1(W, w_0) \to \pi_1(W',w'_0)$ is an isomorphism.
The flow $W' = X' \underset{Y}{\times} Z$ has a fiber bundle structure with $Z$ as a basis
and $\T^{\nu-3} \times K$ as the fiber space.
Thus with $M = {\psi'}^{-1}(z_0)$ we get from the inclusion maps sequence
$$
(M, w'_0) \overset{i}{\to} (W', w'_0) \overset{j}{\to} (W', M, w'_0),
$$ 
the exact sequence
$$
\Z^{\nu - 3}= \pi_1(M, w'_0)  \overset{i_*}{\to} \pi_1(W', w'_0)\overset{j_*}{\to}
\pi_1(W', M, w'_0) \cong \pi_1(Z, z_0).
$$
It follows that $\pi_1(W', w'_0)$, which is abelian, is also a $\Z^{\nu- 3}$-extension of the non-abelian
group $\pi_1(Z, z_0) = \Ga$. This contradiction completes the proof.
\end{proof}

\begin{thm}
The minimal flow $(Z,T)$ is relatively disjoint  from $(|\Wcal|, T)$
over its largest equicontinuous factor $(Y,T)$.
\end{thm}

\begin{proof}
We have seen that the directed collection of flows $\{X(\beta_1, \dots, \beta_u; L_1, \ldots, L_u) \}$
generates $|\Wcal|$ as its inverse limit.
(The assumption $\beta_1 = \al, \ \beta_2 = \beta$ causes no loss in generality.)
The theorem now follows from Lemma \ref{lem-disj} and the fact that relative disjointness 
over a fixed factor is preserved under inverse limits.
\end{proof}

\end{document}